\documentclass[12pt,draft]{article}

\usepackage[cp1251]{inputenc}
\usepackage[T2A]{fontenc}
\usepackage[english,ukrainian]{babel}
\usepackage{amsmath,amsfonts,amssymb}
\usepackage{geometry}

\begin{document}

\begin{center}

\large

\textbf{LAWRUK ELLIPTIC BOUNDARY-VALUE PROBLEMS\\ FOR  HOMOGENEOUS DIFFERENTIAL EQUATIONS}

\medskip

\textbf{A.V. Anop}

\medskip

\normalsize

Institute of Mathematics, National Academy of Sciences of Ukraine, Kyiv

\medskip

E-mail: ahlv@ukr.net

\bigskip\bigskip

\large

\textbf{ЕЛІПТИЧНІ ЗА ЛАВРУКОМ КРАЙОВІ ЗАДАЧІ\\ДЛЯ ОДНОРІДНИХ ДИФЕРЕНЦІАЛЬНИХ РІВНЯНЬ}\footnote{Публікація містить результати досліджень, проведених за грантом Президента України за конкурсним проектом Ф75/29007 Державного фонду фундаментальних досліджень.}

\medskip

\textbf{А.В. Аноп}

\normalsize

\end{center}

\medskip

\noindent We investigate Lawruk elliptic boundary-value problems for homogeneous differential equations in a two-sided refined Sobolev scale. These problems contain additional unknown functions in the boundary conditions of arbitrary orders. The scale consists of inner-product H\"ormander spaces whose orders of regularity are given by any real number and a function which varies slowly at infinity in the sense of Karamata. We establish theorems on the Fredholm property of the problems in the refined Sobolev scale and on local regularity and local \textit{a priori} estimate (up to the boundary of the domain) of their generalized solutions. We find sufficient conditions under which  components of these solutions are $l\geq0$ times continuously differentiable functions.

\medskip

\noindent\textbf{Keywords}: elliptic boundary-value problem, refined Sobolev scale, Fredholm operator, regularity of solution, a priori estimate.

\bigskip

\noindent \textbf{Вступ.} У різних застосуваннях, зокрема, у теорії пружності і гідродинаміці виникають крайові задачі з додатковими невідомими функціями у крайових умовах. Широкий клас таких задач еліптичного типу увів Б. Лаврук [\ref{Lawruk63a}]. На відміну від еліптичних за Я.Б.~Лопатинським крайових задач, цей клас замкнений  відносно переходу до формально спряженої задачі. Еліптичні за Лавруком крайові задачі є нетеровими на відповідних парах просторів Соболєва\,--\,Ройтберга довільного дійсного порядку [\ref{KozlovMazyaRossmann97}, \ref{Roitberg99}]. Такі простори, окрім розподілів в евклідовій області, містять також вектори іншої природи, що ускладнює їх застосування.

Мета цієї роботи --- встановити теореми про нетеровість  еліптичних за Лавруком крайових задач і властивості їх узагальнених розв'язків у  двобічній шкалі просторів, які складаються виключно з розподілів в області. У роботі розглянуто важливий випадок однорідних еліптичних рівнянь. Їх розв'язки досліджуємо в уточненій соболєвській шкалі, яка складається з гільбертових ізотропних просторів Хермандера [\ref{Hermander65}, п.~2.2; \ref{Hermander86}, п.~10.1], для яких показниками регулярності розподілів служать довільні дійсне число і додатна функція, повільно змінна на нескінченності за Й. Караматою. Ця  шкала була виділена і досліджена В.А. Михайлецем і О.О. Мурачем і застосована до дослідження еліптичних за Лопатинським крайових задач [\ref{MikhailetsMurach06UMJ3}, \ref{MikhailetsMurach06UMG11}, \ref{MikhailetsMurach10}, \ref{MikhailetsMurach12BJMA2}]. Вона містить двобічну гільбертову шкалу соболєвських просторів. У цій роботі, на відміну від [\ref{ChepuruhinaDop15}], не накладаємо обмеження на порядки крайових операторів і досліджуємо локальні (аж до межі області) властивості узагальнених розв'язків еліптичних за Лавруком крайових задач.

\textbf{1. Постановка задачі.}  Нехай $\Omega$~--- довільна  обмежена  область у евклідовому просторі
$\mathbb{R}^{n}$, де $n\geq2$, з нескінченно гладкою межею $\Gamma$.
Нехай задано цілі числа $q\geq1$, $\varkappa\geq1$, $m_{1},\ldots,m_{q+\varkappa}$ і $r_{1},\ldots,r_{\varkappa}$.
Розглянемо в $\Omega$  крайову задачу, яка складається з однорідного еліптичного рівняння
\begin{equation}\label{1f1}
Au=0\quad\mbox{в}\quad\Omega
\end{equation}
порядку $2q$ і крайових умов
\begin{equation}\label{1f2}
B_{j}u+\sum_{k=1}^{\varkappa}C_{j,k}v_{k}=g_{j}\quad\mbox{на}\quad\Gamma,
\quad j=1,...,q+\varkappa.
\end{equation}
У ній є невідомими функція $u$ в області $\Omega$ і $\varkappa$ функцій $v_{1},\ldots,v_{\varkappa}$ на $\Gamma$ (всі функції та розподіли вважаємо комплекснозначними, а функціональні простори --- комплексними). Тут $A:=A(x,D)$~--- лінійний диференціальний оператор на $\overline{\Omega}:=\Omega\cup\Gamma$ порядку $2q$, кожне $B_{j}:=B_{j}(x,D)$ є крайовим лінійним диференціальним оператором на $\Gamma$ порядку $\mathrm{ord}\,B_{j}\leq m_{j}$, а кожне $C_{j,k}=C_{j,k}(x,D_{\tau})$ є дотичним лінійним диференціальним оператором на $\Gamma$ порядку $\mathrm{ord}\,C_{j,k}\leq m_{j}+r_{k}$. Усі коефіцієнти цих диференціальних операторів є нескінченно гладкими функціями, заданими на $\overline{\Omega}$ і $\Gamma$ відповідно. (Звісно, диференціальні оператори від'ємних порядків вважаються нуль-операторами).

Припускаємо (і це природно), що
$$
m:=\max \{m_{1},\ldots,m_{q+\varkappa}\}=\max\{\mathrm{ord}B_1,\ldots,\mathrm{ord}B_{q+\varkappa}\}
$$
і $m\geq-r_{k}$  для кожного $k\in\{1,\ldots,\varkappa\}$.
Покладемо $\mu:=\max\{m+1,2q\}$.

Припускаємо, що крайова задача \eqref{1f1}, \eqref{1f2} є еліптичною в області $\Omega$ за Б.~Лавруком [\ref{Lawruk63a}], тобто диференціальний оператор $A$ правильно еліптичний на $\overline{\Omega}$, а система крайових умов \eqref{1f2} задовольняє аналог умови Шапiро\,--\,Лопатинського щодо $A$ на $\Gamma$ (див., наприклад, [\ref{KozlovMazyaRossmann97}, п.~3.1.2]).

Позначимо через $C^\infty(\overline{\Omega},A)$ множину всіх функцій $u\in C^\infty(\overline{\Omega})$ таких, що
$Au=0$ в області $\Omega$.
Пов'яжемо із задачею \eqref{1f1}, \eqref{1f2} лінійне відображення
\begin{equation}\label{1f3}
\begin{gathered}
\mathcal{L}:(u,v_{1},...,v_{\varkappa})\to
\biggl(B_{1}u+\sum_{k=1}^{\varkappa}C_{1,k}v_{k},...,
B_{q+\varkappa}u+\sum_{k=1}^{\varkappa}C_{q+\varkappa,k}v_{k}\biggr),\\
\mbox{де}\quad u\in C^{\infty}(\overline{\Omega},A),\quad
v_{1},\ldots,v_{\varkappa}\in C^{\infty}(\Gamma).
\end{gathered}
\end{equation}
 Досліджуємо властивості його продовження за неперервністю у парах гільбертових функціональних просторів, які утворюють уточнену соболєвську шкалу [\ref{MikhailetsMurach10}, п. 1.3.3].

Для опису області значень цього продовження  потрібна така спеціальна формула Гріна [\ref{KozlovMazyaRossmann97}, формула (4.1.10)]:
\begin{gather*}
(Au,\omega)_{\Omega}+\sum_{j=1}^{\mu-2q}(D_{\nu}^{j-1}Au,w_{j})_{\Gamma}+
\sum_{j=1}^{q+\varkappa}\biggl(B_{j}u+
\sum_{k=1}^{\varkappa}C_{j,k}v_{k},h_{j}\biggr)_{\Gamma}=\\
=(u,A^{+}\omega)_{\Omega}+\sum_{k=1}^{\mu}\biggl(D_{\nu}^{k-1}u,K_{k}\omega+
\sum_{j=1}^{\mu-2q}R_{j,k}^{+}w_{j}+
\sum_{j=1}^{q+\varkappa}Q_{j,k}^{+}h_{j}\biggr)_{\Gamma}+\sum_{k=1}^{\varkappa}\biggl(v_{k},
\sum_{j=1}^{q+\varkappa}C_{j,k}^{+}h_{j}\biggr)_{\Gamma},
\end{gather*}
для довільних функцій $u,\omega\in C^{\infty}(\overline{\Omega})$ і $v_{1},\ldots,v_{\varkappa},w_{1},\ldots,w_{m-2q+1},h_{1},\ldots,h_{q+\varkappa}\in C^{\infty}(\Gamma)$.
Тут і далі $D_\nu:=i\partial/\partial\nu$, де $\nu$ --- поле ортів внутрішніх нормалей до межі $\Gamma$, а через $(\cdot,\cdot)_{\Omega}$ і $(\cdot,\cdot)_{\Gamma}$ позначено відповідно скалярні добутки у гільбертових просторах $L_{2}(\Omega)$ і $L_{2}(\Gamma)$ функцій квадратично інтегровних на $\Omega$ і $\Gamma$ відносно мір Лебега, а також продовженням за неперервністю цих скалярних добутків. Окрім того, $A^{+}$~--- диференціальний оператор, формально спряжений до $A$ відносно $(\cdot,\cdot)_{\Omega}$, усі $C_{j,k}^{+}$, $R_{j,k}^{+}$ і $Q_{j,k}^{+}$ --- дотичні диференціальні оператори, формально спряжені  відповідно до $C_{j,k}$, $R_{j,k}$ і $Q_{j,k}$ відносно $(\cdot,\cdot)_{\Gamma}$. Тут дотичні лінійні диференціальні оператори $R_{j,k}:=R_{j,k}(x,D_{\tau})$ і $Q_{j,k}:=Q_{j,k}(x,D_{\tau})$ узяті із зображення крайових диференціальних операторів $D_{\nu}^{j-1}A$ і $B_{j}$ у вигляді
\begin{equation*}
D_{\nu}^{j-1}A(x,D)=\sum_{k=1}^{m+1}R_{j,k}(x,D_{\tau})D_{\nu}^{k-1},\quad
B_{j}(x,D)=\sum_{k=1}^{m+1}Q_{j,k}(x,D_{\tau})D_{\nu}^{k-1}.
\end{equation*}
Нарешті, кожне $K_{k}:=K_{k}(x,D)$~--- деякий крайовий лінійний диференціальний оператор на $\Gamma$ порядку $\mathrm{ord}\,K_{k}\leq2q-k$ з коефіцієнтами класу $C^{\infty}(\Gamma)$.

Спеціальна формула Гріна приводить до такої крайової
задачі в області $\Omega$:
\begin{equation}\label{1f4}
A^{+}\omega=0\quad\mbox{в}\quad\Omega,
\end{equation}
\begin{equation}\label{1f5}
\begin{gathered}
K_{k}\omega+\sum_{j=1}^{\mu-2q}R_{j,k}^{+}w_{j}+
\sum_{j=1}^{q+\varkappa}Q_{j,k}^{+}h_{j}=\psi_{k}\quad
\mbox{на}\quad\Gamma,\quad k=1,...,\mu,
\end{gathered}
\end{equation}
\begin{equation}\label{1f6}
\sum_{j=1}^{q+\varkappa}C_{j,k}^{+}h_{j}=\psi_{\mu+k}
\quad\mbox{на}\quad\Gamma,\quad k=1,...,\varkappa.
\end{equation}
Ця задача містить, окрім невідомої функції $\omega$ в області $\Omega$, ще $\mu-q+\varkappa$ додаткових невідомих функцій $w_{1},\ldots,w_{\mu-2q}$ і $h_{1},\ldots,h_{q+\varkappa}$ на межі $\Gamma$. Задачу \eqref{1f4}\,--\,\eqref{1f6} називають формально спряженою до задачі \eqref{1f1}, \eqref{1f2} відносно розглянутої спеціальної формули Гріна.  Відомо [\ref{KozlovMazyaRossmann97}, теорема~4.1.1], що еліптичність за Лавруком задачі \eqref{1f1}, \eqref{1f2} рівносильна еліптичності за Лавруком формально спряженої задачі \eqref{1f4}\,--\,\eqref{1f6}.

\textbf{2. Уточнена соболєвська шкала} складається з гільбертових просторів Хермандера, для яких показником регулярності розподілів служить пара параметрів --- числовий $s\in\mathbb{R}$ і функціональний $\varphi\in\mathcal{M}$. Тут $\mathcal {M}$ --- множина всіх вимірних за Борелем функцій $\varphi: [1,\infty)\rightarrow(0,\infty)$, які обмежені і відокремлені від нуля на кожному компакті і повільно змінюються на нескінченності за Й. Караматою, тобто $\varphi(\lambda t)/\varphi(t)\rightarrow 1$ при $t\rightarrow\infty$ для кожного $\lambda>0$ (див., наприклад, [\ref{Seneta85}]). Цю шкалу ввели і дослідили В.А. Михайлець і О.О. Мурач [\ref{MikhailetsMurach06UMJ3}] (див. також їх монографію [\ref{MikhailetsMurach10}] і огляд [\ref{MikhailetsMurach12BJMA2}]).

Нехай $s\in\mathbb{R}$ і $\varphi\in\mathcal {M}$. Означимо простір
$H^{s,\varphi}(\cdot)$ спочатку на $\mathbb {R}^{n}$, а потім на $\Omega$ і $\Gamma$. Будемо дотримуватися монографії [\ref{MikhailetsMurach10}, пп.~1.3, 2.1, 3.2].

За означенням, комплексний лінійний простір $H^{s,\varphi}(\mathbb{R}^{n})$, де $n\in\mathbb{N}$, складається з усіх повільно зростаючих розподілів $w$ на $\mathbb{R}^{n}$ таких, що функція $\langle\xi\rangle^{s}\varphi(\langle\xi\rangle)\widehat{w}(\xi)$ аргументу $\xi\in\mathbb {R}^{n}$ квадратично інтегровна на $\mathbb {R}^{n}$ за Лебегом. Тут $\widehat{w}$ --- перетворення Фур'є розподілу $w$, а $\langle\xi\rangle:=(1+|\xi|^{2})^{1/2}$. У просторі $H^{s,\varphi}(\mathbb {R}^{n})$ означено скалярний добуток  і норму за формулами
$$
(w_{1},w_{2})_{s,\varphi,\mathbb {R}^{n}}:= \int\limits_{\mathbb{R}^{n}}
\langle\xi\rangle^{2s}\varphi^{2}(\langle\xi\rangle)
\widehat{w_{1}}(\xi)\overline{\widehat{w_{2}}(\xi)}d\xi,\quad \|w\|_{s,\varphi,\mathbb {R}^{n}}:=(w,w)_{s,\varphi,\mathbb {R}^{n}}^{1/2}.
$$
Простір $H^{s,\varphi}(\mathbb{R}^{n})$ є ізотропним гільбертовим випадком простору $\mathcal{B}_{p,\mu}$, введеного і дослідженого Л.~Хермандером [\ref{Hermander65}, п.~2.2] (див. також [\ref{Hermander86}, п.~10.1). А саме, $H^{s,\varphi}(\mathbb{R}^{n})=\mathcal{B}_{2,\mu}$, якщо  $\mu(\xi)=\langle\xi\rangle^{s}\varphi(\langle\xi\rangle)$ для довільного $\xi\in\mathbb{R}^{n}$.  Простір $\mathcal{B}_{2,\mu}=H^\mu$ і його версії для евклідових областей ввели і дослідили також Л.Р. Волевич і Б. П. Панеях [\ref{VolevichPaneah65}, \S 2].

Якщо $\varphi(\cdot)\equiv1$, то $H^{s,\varphi}(\mathbb {R}^{n})$ стає гільбертовим простором Соболєва $H^{s}(\mathbb {R}^{n})$ порядку $s$. У загальній ситуації виконуються неперервні та щільні вкладення
$H^{s+\varepsilon}(\mathbb{R}^{n})\hookrightarrow H^{s,\varphi}(\mathbb{R}^{n})\hookrightarrow H^{s-\varepsilon}(\mathbb{R}^{n})$ для довільного $\varepsilon>0$. З огляду на це клас функціональних просторів  $H^{s,\varphi}(\mathbb{R}^{n})$, де $s\in\mathbb{R}$ і $\varphi\in\mathcal {M}$, називають уточненою соболєвською шкалою на $\mathbb {R}^{n}$. Її аналоги для евклідової області $\Omega$ і замкненого
компактного многовиду $\Gamma$ вводяться у стандартний спосіб.
Наведемо відповідні означення.

Лінійний простір $H^{s,\varphi}(\Omega)$ складається зі звужень в
область $\Omega$ усіх розподілів $w\in H^{s,\varphi}(\mathbb{R}^{n})$ і наділений нормою
$$
\|u\|_{s,\varphi,\Omega}:=
\inf\bigl\{\,\|w\|_{s,\varphi,\mathbb{R}^{n}}:\,
w\in
H^{s,\varphi}(\mathbb{R}^{n}),\;w=u\;\,\mbox{в}\;\,\Omega\,\bigr\},
$$
де $u\in H^{s,\varphi}(\Omega)$. Цей простір гільбертів і сепарабельний відносно вказаної норми. Множина $C^{\infty}(\overline{\Omega})$ щільна у ньому.

Лінійний простір $H^{s,\varphi}(\Gamma)$ складається з усіх розподілів на многовиді $\Gamma$, які в локальних координатах дають елементи простору $H^{s,\varphi}(\mathbb{R}^{n-1})$. Дамо детальне означення. Довільно виберемо скінченний атлас із $C^{\infty}$-структури на многовиді $\Gamma$, утворений локальними картами $\pi_j:\mathbb{R}^{n-1}\leftrightarrow\Gamma_{j}$, де $j=1,\ldots,\lambda$.
Тут відкриті множини $\Gamma_{1},\ldots,\Gamma_{\lambda}$ складають скінченне покриття многовиду $\Gamma$. Нехай, окрім того, функції $\chi_j\in C^{\infty}(\Gamma)$, де $j=1,\ldots,\lambda$, утворюють розбиття одиниці на $\Gamma$, що задовольняє умову $\mathrm{supp}\,\chi_j\subset\Gamma_j$. Тоді, за означенням, простір $H^{s,\varphi}(\Gamma)$ складається з усіх розподілів $h$ на $\Gamma$ таких, що $(\chi_{j}h)\circ\pi_{j}\in H^{s,\varphi}(\mathbb{R}^{n-1})$
для кожного номера $j\in\{1,\ldots,\lambda\}$, де  $(\chi_{j}h)\circ\pi_{j}$ є зображенням розподілу $h$ у локальній карті $\pi_{j}$. У просторі $H^{s,\varphi}(\Gamma)$ означено норму за формулою
$$
\|h\|_{s,\varphi,\Gamma}:=\biggl(\,\sum_{j=1}^{\lambda}
\|(\chi_{j}h)\circ\pi_{j}\|_{s,\varphi,\mathbb{R}^{n-1}}^{2}
\biggr)^{1/2}.
$$
Він є гільбертовим і сепарабельним відносно цієї норми, та з точністю до еквівалентності норм не залежить від вибору атласу і розбиття одиниці [\ref{MikhailetsMurach10}, теорема~2.3]. Множина $C^{\infty}(\Gamma)$ щільна у цьому просторі $H^{s,\varphi}(\Gamma)$.

Якщо $\varphi(\cdot)\equiv1$, то $H^{s,\varphi}(\Upsilon)=:H^{s}(\Upsilon)$ є простір Соболєва порядку $s\in\mathbb{R}$ на $\Upsilon\in\{\Omega,\Gamma\}$. У загальній ситуації виконуються компактні та щільні вкладення $H^{s+\varepsilon}(\Upsilon)\hookrightarrow H^{s,\varphi}(\Upsilon)\hookrightarrow H^{s-\varepsilon}(\Upsilon)$ для довільного $\varepsilon>0$

\textbf{3. Основні результати} стосуються властивостей еліптичної крайової задачі \eqref{1f1}, \eqref{1f2} в уточненій соболєвській шкалі.
Позначимо через $N$ лінійний простір усіх розв'язків $(u,v_{1},\ldots,v_{\varkappa})\in C^{\infty}(\overline{\Omega})\times (C^{\infty}(\Gamma))^{\varkappa}$ цієї задачі у випадку, коли  кожне $g_{j}=0$ на~$\Gamma$. Аналогічно, позначимо через $N^{+}$ лінійний простір усіх розв'язків
$$
(\omega,w_{1},\ldots,w_{\mu-2q},h_{1},\ldots,h_{q+\varkappa})\in C^{\infty}(\overline{\Omega})\times(C^{\infty}(\Gamma))^{\mu-2q}\times(C^{\infty}(\Gamma))^{q+\varkappa}
$$
формально спряженої крайової задачі \eqref{1f4}\,--\,\eqref{1f6} у випадку, коли всі $\psi_{k}=0$ і $\psi_{m+k}=0$ на~$\Gamma$. Також позначимо через $N_{1}^{+}$ простір усіх векторів $(h_{1},\ldots,h_{q+\varkappa})\in (C^{\infty}(\Gamma))^{q+\varkappa}$,
для яких існують функції $\omega\in C^{\infty}(\overline{\Omega})$ і $w_{1},\ldots,w_{\mu-2q}\in C^{\infty}(\Gamma)$ такі, що $(\omega,w_{1},\ldots,w_{\mu-2q},h_{1},\ldots,h_{q+\varkappa})\in N^+$. Оскільки обидві задачі еліптичні за Лавруком в $\Omega$, то простори $N$, $N^{+}$ та $N_{1}^{+}$ скінченновимірні [\ref{KozlovMazyaRossmann97}, наслідок~4.1.1].

Для довільних $s\in\mathbb{R}$ і $\varphi\in\mathcal{M}$ позначимо через $H^{s,\varphi}(\Omega,A)$ множину усіх розподілів $u\in H^{s,\varphi}(\Omega)$ таких, що $Au=0$ в $\Omega$. Розглядаємо $H^{s,\varphi}(\Omega,A)$ як (замкнений) підпростір гільбертового простору $H^{s,\varphi}(\Omega)$. Множина  $C^{\infty}(\overline{\Omega},A)$ щільна в $H^{s,\varphi}(\Omega,A)$ згідно з [\ref{MikhailetsMurach10}, теорема 3.11].

\textbf{Теорема 1.} \it Для довільних $s\in\mathbb{R}$ і $\varphi\in\mathcal{M}$ відображення \eqref{1f3} продовжується єдиним чином (за неперервністю) до обмеженого оператора
\begin{equation}\label{4f12}
\begin{gathered}
\mathcal{L}:\,H^{s,\varphi}(\Omega,A)\oplus
\bigoplus_{k=1}^{\varkappa}H^{s+r_{k}-1/2,\varphi}(\Gamma)\rightarrow \bigoplus_{j=1}^{q+\varkappa}H^{s-m_{j}-1/2,\varphi}(\Gamma).
\end{gathered}
\end{equation}
Цей оператор нетерів. Його ядро збігається з $N$, а область значень складається з усіх векторів $(g_{1},\ldots,g_{q+\varkappa})
\in\bigoplus_{j=1}^{q+\varkappa}H^{s-m_{j}-1/2,\varphi}(\Gamma)$ таких, що
\begin{equation}\label{4f12c}
(g_{1},h_{1})_{\Gamma}+\cdots+(g_{q+\varkappa},h_{q+\varkappa})_{\Gamma}=0\quad\mbox{для кожного}\quad
(h_{1},\ldots,h_{q+\varkappa})\in N^{+}_{1}.
\end{equation}
Індекс оператора \eqref {4f12} дорівнює $\mathrm{dim}N - \mathrm{dim}N_{1}^{+}$ і не залежить від $s$ та~$\varphi$. \rm

Нагадаємо, що лінійний обмежений оператор $T:E_{1}\rightarrow E_{2}$, який діє у парі банахових просторів $E_{1}$ і $E_{2}$, називають нетеровим, якщо його ядро $\ker T$ і коядро $E_{2}/T(E_{1})$ скінченновимірні. Нетерів оператор має замкнену область значень $T(E_{1})$  та скінченний індекс $\mathrm{ind}\,T:=\dim\ker T-\dim(E_{2}/T(E_{1}))$.

Перейдемо до локальних властивостей узагальнених розв'язків досліджуваної задачі \eqref{1f1}, \eqref{1f2}. Наведемо спочатку їх означення. Позначимо через $S'(\Omega)$ лінійний топологічний простір звужень в область $\Omega$ усіх повільно зростаючих розподілів на $\mathbb{R}^n$ і покладемо $S'(\Omega,A):=\{u\in S'(\Omega):Au=0\;\mbox{в}\;\Omega\}$. Позначимо також через $D'(\Gamma)$ лінійний топологічний простір усіх розподілів на $\Gamma$. Нехай
\begin{equation}\label{4f18}
(u,v):=(u,v_{1},\ldots,v_{\varkappa})\in
S'(\Omega,A)\times(D'(\Gamma))^\varkappa.
\end{equation}
Вектор $(u,v)$ називаємо узагальненим розв'язком цієї задачі, де $g=(g_1,\ldots,g_{q+\varkappa})\in(D'(\Gamma))^{q+\varkappa}$, якщо $\mathcal{L}(u,v)=g$, де $\mathcal{L}$ --- оператор \eqref{4f12} для деяких достатньо малого $s$ і $\varphi$. У цьому випадку $u\in C^{\infty}(\Omega)$, оскільки рівняння \eqref{1f1} еліптичне в $\Omega$. У термінах уточненої соболєвської шкали дамо достатні умови регулярності розв'язку $u$ в області $\Omega$ впритул до куска $\Gamma_0$ її межі.

Нехай $\Gamma_{0}$ --- довільна непорожня підмножина многовиду $\Gamma$. Позначимо через $H^{\sigma,\varphi}_{\mathrm{loc}}(\Omega,\Gamma_{0})$, де $\sigma\in\mathbb{R}$ і $\varphi\in\mathcal{M}$, лінійний простір усіх розподілів $u\in S'(\Omega)$ таких, що $\chi u\in H^{\sigma,\varphi}(\Omega)$ для довільної функції $\chi\in C^{\infty}(\overline{\Omega})$, носій якої   лежить в  $\Omega\cup\Gamma_{0}$. Аналогічно, позначимо через $H^{\sigma,\varphi}_{\mathrm{loc}}(\Gamma_{0})$ лінійний простір усіх розподілів $h\in D'(\Gamma)$ таких, що $\chi_1 h\in H^{\sigma,\varphi}(\Gamma)$ для довільної функції $\chi_1\in C^{\infty}(\Gamma)$, носій якої   лежить в  $\Gamma_{0}$.

\textbf{Теорема 2.} \it
Нехай $s\in\mathbb{R}$ і $\varphi\in\mathcal{M}$, а вектор \eqref{4f18} є узагальненим розв'язком еліптичної крайової задачі \eqref{1f1}, \eqref{1f2}, праві частини якої задовольняють умову $g_j\in H^{s-m_j-1/2,\varphi}_\mathrm{loc}(\Gamma_{0})$ для кожного $j\in \{1,\ldots,q+\varkappa\}$. Тоді $u\in H^{s,\varphi}_{\mathrm{loc}}(\Omega,\Gamma_{0})$ і $v_k\in H^{s+r_k-1/2,\varphi}_\mathrm{loc}(\Gamma_{0})$ для кожного $k\in \{1,\ldots,\varkappa\}$.
\rm

Цей розв'язок задовольняє таку апріорну оцінку:

\textbf{Теорема 3.} \it
Нехай число $\lambda>0$, а функції  $\chi,\,\eta\in C^\infty(\overline{\Omega})$ такі, що  їх носії лежать в $\Omega\cup\Gamma_0$ і $\eta=1$ в околі $\mathrm{supp}\,\chi$; покладемо $\chi_1:=\chi\!\upharpoonright\!\Gamma$ і $\eta_1:=\eta\!\upharpoonright\!\Gamma$. Тоді існує число $c>0$ таке, що
\begin{equation}\label{1f19}
\begin{gathered}
\|\chi u\|_{s,\varphi,\Omega}+\sum_{k=1}^{\varkappa}\|\chi_1 v_k\|_{s+r_k-1/2,\varphi,\Gamma}\leq c\,\sum_{j=1}^{q+\varkappa}\|\eta_1g_j\|_{s-m_j-1/2,\varphi,\Gamma}+\\
+c\,\Bigl(\,\|\eta u\|_{s-\lambda,\varphi,\Omega}
+\sum_{k=1}^{\varkappa}\|\eta_1v_k\|_{s+r_k-1/2-\lambda,\varphi,\Gamma}\Bigl)
\end{gathered}
\end{equation}
для довільних векторів $\eqref{4f18}$ і $(g_1,\ldots,g_{q+\varkappa})$ з теореми $2$.
\rm

До цих теорем є найбільш близькими результати робіт [\ref{MikhailetsMurach06UMG11}, \ref{ChepuruhinaDop15}, \ref{AnopMurachDop18}], у яких досліджено еліптичні крайові задачі для однорідного рівняння \eqref{1f1} у шкалах гільбертових просторів Хермандера. В [\ref{MikhailetsMurach06UMG11}] (див. також монографію [\ref{MikhailetsMurach10}, п. 3.3]) досліджено регулярну еліптичну крайову задачу, а в [\ref{ChepuruhinaDop15}]~--- еліптичну за Лавруком задачу \eqref{1f1}, \eqref{1f2}  за додаткового припущення, що $\mathrm{ord}\,B_j\leq2q-1$. Для цих задач доведено теорему 1, а також теореми 2 і 3 у випадку, коли $\Gamma_0=\Gamma$ (глобальна регулярність розв'язку) і $\chi=\eta=1$ на $\overline{\Omega}$ (глобальна апріорна оцінка розв'язку). В [\ref{AnopMurachDop18}] встановлено версії теорем 1 і 2 для еліптичних за Лопатинським крайових задач (без додаткових невідомих функцій у крайових умовах), які розглядаються в розширеній соболєвській шкалі. Якщо еліптичне рівняння \eqref{1f1} неоднорідне, то аналоги теорем 1--3 правильні за умови, що $s>m+1/2$ [\ref{ChepuruhinaKasirenko17}].

\textbf{4. Застосування.} Одним із застосувань уточненої соболєвської шкали є достатні умови, за яких компоненти  узагальненого розв'язку \eqref{4f18} еліптичної задачі \eqref{1f1}, \eqref{1f2} належать до просторів $l$ разів неперервно диференційовних функцій на $\Omega\cup\Gamma_{0}$ і $\Gamma_0$ відповідно. Ці умови випливають із теореми 2 і такої версії теореми вкладення Хермандера [\ref{Hermander65}, теорема 2.2.7]: кожне з вкладень  $H^{l+n/2,\varphi}(\Omega)\subset C^{l}(\overline{\Omega})$ і
$H^{l+(n-1)/2,\varphi}(\Gamma)\subset C^{l}(\Gamma)$, де $0\leq l\in\mathbb{Z}$ і $\varphi\in\mathcal{M}$, еквівалентне умові
\begin{equation}\label{1f45}
\int_{1}^{\infty}\frac{dt}{t\,\varphi^{2}(t)}<\infty;
\end{equation}
ці вкладення компактні (див. [\ref{MikhailetsMurach10}, теореми 2.8 і 3.4]).

\textbf{Теорема 4.} \it Нехай цілі числа $l\geq0$ і $k\in\{1,...,\varkappa\}$. Якщо вектор \eqref{4f18}  задовольняє умову теореми $2$, де $s=l+n/2$, а $\varphi\in\mathcal{M}$ --- умову \eqref{1f45}, то $u\in C^{l}(\Omega\cup\Gamma_{0})$. Якщо вектор \eqref{4f18}  задовольняє умову теореми $2$, де $s=l-r_k+n/2$, а $\varphi\in\mathcal{M}$ --- умову \eqref{1f45}, то $v_{k}\in C^{l}(\Gamma_0)$.
\rm

Важливо, що умова \eqref{1f45} є точною у цій теоремах.

Наведемо також застосування отриманих результатів до деякого узагальнення задачі Пуанкаре\,--\,Стєклова [\ref{QuarteroniValli91}]. Окрім задачі \eqref{1f1}, \eqref{1f2}, розглянемо еліптичну за Лавруком задачу, яка складається з того ж самого однорідного диференціального рівняння \eqref{1f1} і крайових умов
\begin{equation}\label{1f2b}
\widetilde{B}_{j}u+\sum_{k=1}^{\varkappa}\widetilde{C}_{j,k}v_{k}=\widetilde{g}_{j}\quad\mbox{на}\quad\Gamma,
\quad j=1,...,q+\varkappa.
\end{equation}
Тут, подібно до \eqref{1f2}, кожне $\widetilde{B}_{j}$ є крайовим лінійним диференціальним оператором на $\Gamma$ порядку $\mathrm{ord}\,\widetilde{B}_{j}\leq \widetilde{m}_{j}$, а кожне $\widetilde{C}_{j,k}$ є дотичним лінійним диференціальним оператором на $\Gamma$ порядку $\mathrm{ord}\,\widetilde{C}_{j,k}\leq \widetilde{m}_{j}+r_{k}$, де $\widetilde{m}_{1},\ldots,\widetilde{m}_{q+\varkappa}$ --- цілі числа.
Покладемо $g:=(g_1,\ldots,g_{q+\varkappa})$ і $\widetilde{g}=(\widetilde{g}_1,\ldots,\widetilde{g}_{q+\varkappa})$.

Повя'жемо з цими задачами відображення $\Psi$, яке вектору $g\in(D'(\Gamma))^{q+\varkappa}$, підпорядкованому умові \eqref{4f12c}, ставить у відповідність вектор $\widetilde{g}\in(D'(\Gamma))^{q+\varkappa}$ такий, що $\widetilde{g}$ задовольняє крайові умови \eqref{1f2b}, де $(u,v)$ є узагальненим розв'язком \eqref{4f18} крайової задачі \eqref{1f1}, \eqref{1f2}, підпорядкованим умові $(u,u^{\circ})_\Omega+ (v_{1},v_{1}^{\circ})_{\Gamma}+\ldots+(v_{\varkappa},v_{\varkappa}^{\circ})_{\Gamma}=0$ для кожного
$(u^{\circ},v_{1}^{\circ},\ldots,v_{\varkappa}^{\circ})\in N$. За теоремою 1, це відображення однозначне. Воно є природним узагальненням оператора задачі Пуанкаре\,--\,Стєклова на випадок, коли кількість функцій заданих на межі області може перевищувати половину порядку еліптичного рівняння \eqref{1f1}. Зокрема, якщо $q=1$, $\varkappa=0$, а  \eqref{1f2} і \eqref{1f2b} --- крайові умови Діріхле і Неймана, відображення $\Psi$ є оператором Діріхле\,--\,Неймана.

З теореми 1 випливає, що для будь-яких $s\in\mathbb{R}$ і $\varphi\in\mathcal{M}$ відображення $\Psi:g\mapsto\widetilde{g}$ встановлює ізоморфізм підпростору усіх векторів $g\in\bigoplus_{j=1}^{q+\varkappa}H^{s-m_{j}-1/2,\varphi}(\Gamma)$, які задовольняють умову \eqref{4f12c}, на підпростір усіх векторів $\widetilde{g}\in\bigoplus_{j=1}^{q+\varkappa}H^{s-\widetilde{m}_{j}-1/2,\varphi}(\Gamma)$, які задовольняють умову $(\widetilde{g}_{1},\widetilde{h}_{1})_{\Gamma}+\ldots+(\widetilde{g}_{q+\varkappa},
\widetilde{h}_{q+\varkappa})_{\Gamma}=0$ для кожного $(\widetilde{h}_{1},\ldots,\widetilde{h}_{q+\varkappa})\in \widetilde{N}^{+}_{1}$. Тут $\widetilde{N}^{+}_{1}$ --- аналог скінченновимірного простору $N^{+}_{1}$ для задачі \eqref{1f1}, \eqref{1f2b}.

\bigskip\bigskip

\noindent REFERENCES

\begin{enumerate}

\item\label{Lawruk63a}
Lawruk, B. (1963). Parametric boundary-value problems for elliptic systems of linear differential equations. I. Construction of conjugate problems. Bull. Acad. Polon. Sci. S\'{e}r. Sci. Math.
Astronom. Phys., 11, No. 5, pp. 257–267 (in Russian).

\item\label{KozlovMazyaRossmann97}
Kozlov, V.A., Maz'ya, V.G. \& Rossmann, J. (1997). Elliptic Boundary Value Problems in Domains with Point Singularities. Providence: Amer. Math. Soc.

\item\label{Roitberg99}
Roitberg, Ya.A. (1999). Elliptic boundary value problems in the spaces of distributions. Dordrecht: Kluwer Acad. Publishers.

\item\label{Hermander65}
H\"ormander, L. (1963). Linear Partial Differential Operators. Berlin: Springer-Verlag.

\item\label{Hermander86}
H\"ormander, L. (1983). The Analysis of Linear Partial Differential Operators, vol.~II, Differential Operators with Constant Coefficients. Berlin: Springer-Verlag.

\item\label{MikhailetsMurach06UMJ3}
Mikhailets, V.A. \& Murach, A.A. (2006). Refined scales of spaces and elliptic boundary-value problems. II. Ukr. Math. J., 58, No.~3, pp. 398-417.

\item\label{MikhailetsMurach06UMG11}
Mikhailets, V.A. \& Murach, A.A. (2006). Regular elliptic boundary-value problem for homogeneous equation in two-sided refined scale of spaces. Ukr. Math. J., 58, No.~11, pp. 1748-1767.

\item\label{MikhailetsMurach10}
Mikhailets, V.A. \& Murach A.A. (2014). H\"ormander Spaces, Interpolation, and Elliptic Problems. Berlin/Boston: De Gruyter.

\item\label{MikhailetsMurach12BJMA2}
Mikhailets, V.A. \& Murach, A.A. (2012). The refined Sobolev scale, interpolation, and elliptic problems. Banach J. Math. Anal., 6, No.~2., pp. 211-281.

\item\label{ChepuruhinaDop15}
Chepurukhina, I.S. (2015) A semihomogeneous elliptic problem with additional unknown functions in boundary conditions. Dopov. Nac. akad. nauk. Ukr., No.~7, pp. 20-28 (in Russian).

\item\label{Seneta85}
Seneta, E. (1976). Regularly Varying Functions. Berlin: Springer.

\item\label{VolevichPaneah65}
Volevich, L.R. \& Paneah B.P. (1965). Certain spaces of generalized functions and embedding theorems. Russian Math. Surveys, 20, No.~1, pp. 1-73.

\item\label{ChepuruhinaKasirenko17}
Kasirenko, T.M. \& Chepurukhina, I.S. (2017). Elliptic problems in the sense of Lawruk with boundary operators of higher orders in refined Sobolev scale. Transactions of Institute of Mathematics of NAS of Ukraine, 14, No. 3, pp. 161-204 (in Ukrainian). Kyiv: Institute of Mathematics, NAS of Ukraine.

\item\label{AnopMurachDop18}
Anop, A.V. \& Murach, A.A. (2018) Homogeneous elliptic equations in an extended Sobolev scale. Dopov. Nac. akad. nauk Ukr., No 3, pp. 3-11 (in Ukrainian).

\item\label{QuarteroniValli91}
Quarteroni, A. \& Valli, A. (1991). Theory and application of Steklov\,--\,Poincar\`{e} operators for boundary-value problems. Applied and Industrial Mathematics. Mathematics and Its Applications, Vol. 56, pp. 179--203. Dordrecht: Springer.

\end{enumerate}

\end{document}